\documentclass[12pt]{article}

\usepackage{epsfig}
\usepackage{amssymb}
\usepackage{amsmath}
\usepackage{amsfonts}

\def\diam{{\rm diam}}
\def\dist{{\rm dist}}

\newtheorem{theorem}{Theorem}[section]
\newtheorem{lemma}[theorem]{Lemma}
\newtheorem{corollary}[theorem]{Corollary}

\newtheorem{definition}[theorem]{Definition}

\title{Universal Behavior of Connectivity Properties in Fractal Percolation Models}
\author{
{Erik I. Broman}
\thanks{Department of Mathematics, Chalmers University of Technology,
S-412 96 G\"oteborg, Sweden. E-mail: broman\,@\,chalmers.se}
\,\,\,\,\,\,\,\,\,\,
{Federico Camia}
\thanks{Department of Mathematics, Vrije Universiteit Amsterdam,
De Boelelaan 1081a, 1081 HV Amsterdam, The Netherlands. E-mail:
fede\,@\,few.vu.nl}
}

\begin{document}

\date{}

\maketitle

\begin{abstract}
Partially motivated by the desire to better understand the connectivity
phase transition in fractal percolation, we introduce and study a class
of continuum fractal percolation models in dimension $d \geq 2$. These
include a scale invariant version of the classical (Poisson) Boolean model
of stochastic geometry and (for $d=2$) the Brownian loop soup introduced
by Lawler and Werner.

The models lead to random fractal sets whose connectivity properties
depend on a parameter $\lambda$. In this paper we mainly study the
transition between a phase where the random fractal sets are totally
disconnected and a phase where they contain connected components
larger than one point. In particular, we show that there are connected
components larger than one point {\it at} the unique value of $\lambda$
that separates the two phases (called the critical point). We prove that
such a behavior occurs also in Mandelbrot's fractal percolation in all
dimensions $d \geq 2$. Our results show that it is a generic feature,
independent of the dimension or the precise definition of the model,
and is essentially a consequence of scale invariance alone.

Furthermore, for $d=2$ we prove that the presence of connected components
larger than one point implies the presence of a unique, unbounded,
connected component.
\end{abstract}



\medskip \noindent
Submitted to EJP March 23 2009, final version accepted August 31 2010.

\section{Introduction} \label{sec-intro}

Many deterministic constructions generating fractal sets have random analogues
that produce \emph{random fractals} which do not have the self-similarity of
their non-random counterpart, but are \emph{statistically self-similar} in the
sense that enlargements of small parts have the same statistical distribution
as the whole set.

Random fractals can have complex topological structure, for example they
can be highly multiply connected, and can exhibit \emph{connectivity phase
transitions}, corresponding to sudden changes of topological structure as
a continuously varying parameter goes through a critical value.

In this paper, we introduce and study a natural class of random fractals
that exhibit, in dimension $d \geq 2$, such a connectivity phase transition:
when a parameter increases continuously through a critical value, the
connectivity suddenly breaks down and the random fractals become totally
disconnected with probability one. (We remind the reader that a set is
called totally disconnected if it contains no connected component larger
than one point.)
The fractals we study are defined as the complement of the union of sets
generated by a Poisson point process of intensity $\lambda$ times a
\emph{scale invariant} measure on a space of subsets of ${\mathbb R}^d$
(see Section~\ref{sec-main}).

Examples of such random fractals include a scale invariant version
of the classical (Poisson) Boolean model of stochastic geometry (see
\cite{sw,skm}
and \cite{mr} for a multiscale version of the model), the Brownian loop
soup~\cite{lw} (both will be discussed in some more detail in the next
section),
and the models studied in \cite{nacu-werner}.
The scale invariant (Poisson) Boolean model is a natural model for a porous
medium with cavities on many different scales (but it has also been used as a simplified
model in cosmology --- see \cite{fk}). It is obtained via a Poisson point process
in $(d+1)$-dimensional space, where the first $d$ coordinates of the points
give the locations of the centers of $d$-dimensional balls whose radii are
given by the last coordinate. The distribution of the radii $r$ has density
$(1/r)^{d+1}$, which ensures scale invariance.
There is no reason, except simplicity, for using balls, and the model can
be naturally generalized by associating random shapes to the points of the
Poisson process. Another natural way to generalize the model is
obtained by considering a Poisson point process directly in the space of
``shapes,'' i.e., subsets of ${\mathbb R}^d$. In dimension $d=2$, this is how
the Brownian loop soup is defined, with the distribution of the random shapes
given by the distribution of Brownian loops. In this paper we consider this
type of models with general scale invariant distributions on shapes (see
Definition~\ref{def-scale-inv}).
The reason is that we want
to study what features in the behavior of fractal percolation models
are a consequence of scale invariance alone.

Our main result consists in showing that, when the intensity $\lambda$
of the Poisson process is at its critical value, the random fractals are in
the connected phase in the sense that they contain connected components
larger than one point with probability one (see Theorem~\ref{main-thm}).
This is reminiscent of the nature of the phase transition in Mandelbrot's
fractal percolation~\cite{mandelbrot,mandelbrot-book}, which is discussed
in more detail in Section~\ref{sec-mandelbrot}.

Our proof of Theorem~\ref{main-thm} is interesting in that it shows that
the nature of the connectivity phase transition described in the theorem
is essentially a consequence of scale invariance alone, and in particular
does not depend on the dimension $d$. The same proof applies to other
models as well, including Mandelbrot's fractal percolation. (We note that
the proofs of Theorems 5.3 and 5.4 of~\cite{dm} also show the importance
of scale invariance, but are very two-dimensional).
Along the way, we prove a discontinuity result for the probability that a
random fractal contains a connected component crossing a ``shell-like domain,''
which is interesting in its own right (see Corollary~\ref{cor}).

The main result is stated in Section \ref{sec-main} while Sections
\ref{sec-2d} and \ref{sec-mandelbrot} contain additional two-dimensional
results and results concerning Mandelbrot's fractal percolation model
respectively.

\subsection{Two Motivating Examples} \label{sec-examples}
Two prototypical examples of the type of models that we consider in
this paper are a fully scale invariant version of the multiscale
Poisson Boolean model studied in Chapter~8 of~\cite{mr} and in~\cite{mpv1,mpv3}
(see~\cite{beg} for recent results on that model, whose precise definition
is given below) and, in two dimensions, the Brownian loop soup of Lawler
and Werner~\cite{lw}.

The Brownian loop soup with density $\lambda>0$ is a realization
of a Poisson point process with intensity $\lambda$ times the
\emph{Brownian loop measure}, where the latter is essentially the
only measure on loops that is conformally invariant (see~\cite{lw}
and~\cite{lawler-book} for precise definitions). A sample of
the Brownian loop soup is a countable family of unrooted Brownian
loops in ${\mathbb R}^2$ (there is no non-intersection condition
or other interaction between the loops). The Brownian loop measure
can also be considered as a measure on hulls (i.e., compact connected
sets $K \subset {\mathbb R}^2$ such that ${\mathbb R}^2 \setminus K$
is connected) by filling in the bounded loops.

The scale invariant Boolean model in $d$ dimensions is a
Poisson point process on ${\mathbb R}^d \times (0,\infty)$ with
intensity $\lambda \, r^{-(d+1)} \, {\rm d}r \, {\rm d}x$, where $\lambda \in (0,\infty)$,
${\rm d}r$ is Lebesgue measure on ${\mathbb R}^+$ and
${\rm d}x$ is the $d$-dimensional Lebesgue measure. Each realization
$\cal P$ of the point process gives rise to a collection of balls
in ${\mathbb R}^d$ in the following way. For each point
$\xi \in {\cal P}$ there is a corresponding ball $b(\xi)$.
The projection on ${\mathbb R}^d$ of $\xi$ gives the
position of the center of the ball and the radius 
of the ball is given by the value of the last coordinate of $\xi$.

Since we want to show the analogy between the two models, and later
generalize them, we give an alternative description of the scale
invariant Boolean model. One can obtain the random collection
of balls described above as a realization of a Poisson point process with
intensity $\lambda\mu^{Bool}$, where $\mu^{Bool}$ is the measure defined
by $\mu^{Bool}(\tilde E) = \int_D \int_a^b r^{-(d+1)} {\rm d}r {\rm d}x$,
for all sets $\tilde E$ that are collections of balls of radius $r \in (a,b)$
with center in an open subset $D$ of ${\mathbb R}^d$. (Denoting by
$\tilde{\cal E}$ the collection of sets $\tilde E$ used in the definition
of $\mu^{Bool}$, it is easy to see that $\tilde{\cal E}$ is closed under
pairwise intersections. Therefore, the probabilities of events in $\tilde{\cal E}$
determine $\mu^{Bool}$ uniquely as a measure on the $\sigma$-algebra $\sigma(\tilde{\cal E})$.
This choice of $\sigma$-algebra is only an example, later we will work with
different measurable sets.)
Here, $\mu^{Bool}$ plays the same role as the Brownian loop measure in the
definition of the Brownian loop soup.
Note that $\mu^{Bool}$ is scale invariant in the following sense. Let $\tilde E'$
denote the collection of balls with center in $sD$ and radius $r \in (sa,sb)$,
for some scale factor $s$. Then
$\mu^{Bool}(\tilde E') = \int_{sD} \int_{sa}^{sb} r^{-(d+1)} {\rm d}r {\rm d}x
= \int_D \int_a^b (sr)^{-(d+1)} s {\rm d}r s^d {\rm d}x = \mu(\tilde E)$.

We are interested in fractal sets obtained by considering the complement
of the union of random sets like those produced by the scale invariant
Boolean model or the Brownian loop soup,
possibly with a cutoff on the maximal diameter of the random sets.
Fractals have frequently been used to model physical systems, such as porous
media, and in that context the presence of a cutoff is a very
natural assumption.
Furthermore, it will be easy to see from the definitions that
without a cutoff or some other restriction, the complement of the union of the random
sets is a.s.\ empty. An alternative possibility is to consider the
restriction of the scale invariant Boolean model or the Brownian
loop soup to a bounded domain $D$. By this we mean that one keeps
only those balls or Brownian loops that are contained in $D$, which
automatically provides a cutoff on the size of the retained sets.
This approach is particularly natural in the Brownian loop soup
context, since then, when $\lambda$ is below a critical value,
the boundaries of clusters of filled Brownian loops form a realization
of a Conformal Loop Ensemble (see \cite{werner1, werner2, shefwern1}
and \cite{shefwern2} for the definition and properties of Conformal Loop Ensembles).

All proofs are contained in Section \ref{sec-proofs}.

\section{Definitions and Main Results} \label{sec-main}

We first remind the reader of the definition of Poisson point process.
Let $(M,{\cal M},\mu)$ be a measure space, with $M$ a topological space,
$\cal M$ the Borel $\sigma$-algebra, and $\mu$ a $\sigma$-finite
measure. A \emph{Poisson point process with intensity measure $\mu$}
is a collection of random variables $\{N(E): E\in {\cal M}, \ \mu(E)<\infty\}$
satisfying the following properties:
\begin{itemize}
\item With probability 1, $E \mapsto N(E)$ is a counting measure (i.e.,
it takes only nonegative integer values).
\item For fixed $E$, $N(E)$ is a Poisson random variable with mean $\mu(E)$.
\item If $E_1,E_2,\ldots,E_n$ are mutually disjoint, then $N(E_1),N(E_2),\ldots,N(E_n)$
are independent.
\end{itemize}
The random set of points ${\cal P} = \{\xi \in M : N(\{\xi\})=1 \}$
is called a \emph{Poisson realization} of the measure $\mu$.

In the rest of the paper, $d\geq2$, and $D$ will always denote a bounded,
open subset of ${\mathbb R}^d$, which will be called a \emph{domain},
and $\overline{D}$ will denote the closure of $D$.
If $K$ is a subset of ${\mathbb R}^d$, we let
$sK = \{ x \in {\mathbb R}^d : x/s \in K \}$.
Here $M$ will be the set of connected, compact subsets
of ${\mathbb R}^d$ with nonempty interior.
For our purposes, we need not specify the topology, but we
require that the Borel $\sigma$-algebra contains all sets of the form
$E(B;a,b) = \{ K \in M : a < \text{diam}(K) \leq b, K \subset B \}$
for all $0 \leq a < b$ and all Borel sets $B \subset {\mathbb R}^d$.
If we denote the collection of sets
$E(B;a,b)$ by $\cal E$, the latter is a $\pi$-system (i.e., closed under
finite intersections), and one may set ${\cal M}=\sigma(\cal E)$.

We now give a precise definition of scale invariance,
followed by the main definitions of the paper.
\begin{definition} \label{def-scale-inv}
We say that an infinite measure $\mu$ on $(M,{\cal M})$
is \emph{scale invariant} if, for any
$E \in {\cal M}$ with $\mu(E)<\infty$ and any $0<s<\infty$,
$\mu(E')=\mu(E)$, where $E' = \{ K : K/s \in E \}$.
\end{definition}

\begin{definition} \label{def-soup}
A \emph{scale invariant (Poissonian) soup} in $D$ with intensity $\lambda\mu$
is the collection of sets from a Poisson realization of $\lambda\mu$
that are contained in $D$, where $\mu$ is a translation and scale
invariant measure.
\end{definition}
Note that the soup inherits the scale invariance of the measure $\mu$,
so that soup realizations in domains related by uniform scaling are
statistically self-similar. For instance, if $0<s<1$ and $D$ are such
that $sD = \{ x \in {\mathbb R}^d : x/s \in D \} \subset D$, and
${\cal K}_D$ denotes a realization of a scale invariant soup in $D$,
then the collection of sets from ${\cal K}_D$ contained in
$sD$ is distributed like a scaled version $s{\cal K}_{D}$ of ${\cal K}_D$
(where the elements of $s{\cal K}_{D}$ are the sets
$sK = \{ x \in {\mathbb R}^d : x/s \in K \}$ with $K \in {\cal K}_D$).

\begin{definition} \label{def-full-space-soup}
A \emph{full space (Poissonian) soup} with intensity $\lambda\mu$ and cutoff
$\delta>0$ is a Poisson realization from a measure $\lambda\mu_{\delta}$,
where $\mu_{\delta}$ is
the measure induced by $\mu$ on sets of diameter at most
$\delta$, and $\mu$ is a translation and scale invariant measure.
\end{definition}
The scale invariance of the soup can now be expressed in the following
way. Let $0<s<1$, and let $D$ and $D'$ be two disjoint domains such that
$s^{-1} \text{\diam}(D) = \text{\diam}(D') \leq \delta$
(where $\text{diam}(\cdot)$ denotes Euclidean diameter) and with $D'$
obtained by translating $sD$. Then, as before, the sets that are contained
in $D'$ are distributed like a copy scaled by $s$ of the sets contained
in $D$. In other words, the soup is statistically self-similar at all scales
smaller than the cutoff. Note that the full space soup is also stationary
due to the translation invariance of $\mu$.

Clearly, the value $\delta$ of the cutoff is not important, since one can
always scale space to make it become 1. In the rest of the paper, when we
talk about the full space soup without specifying the cutoff $\delta$, we
implicitly assume that $\delta=1$.

We will consider translation and scale invariant measures $\mu$ that satisfy
the following condition.
\begin{itemize}
\item[$(\star)$] Given a domain $D$ and two positive real numbers $d_1<d_2$,
let
$F=F(D;d_1,d_2)$ be the collection of compact connected sets with nonempty
interior that intersect $D$ and have diameters $> d_1$ and $\leq d_2$; then
$\mu(F)<\infty$.
\end{itemize}

\noindent{\bf Remarks} \,
Condition ($\star$) is very natural and is clearly satisfied by $\mu^{Bool}$ and
by the Brownian loop measure (which are also translation and scale invariant).
Its purpose is to ensure that $\lambda_c>0$ in Theorem~\ref{main-thm} below
and to ensure the left-continuity of certain crossing probabilities (see
the beginning of Section~\ref{sec-proofs}).
Note that the set $F$ can be written as $F(D;d_1,d_2) = E(D';d_1,d_2) \cap
E(D' \setminus D;d_1,d_2)^c$, where $D'$ is the (Euclidean) $d_2$-neighborhood
of $D$ and the superscript $c$ denotes the complement. Therefore, $F$ is
measurable by our assumptions on $\cal M$.
\vspace{3mm}

We are now ready to state the main results of the paper.

\begin{theorem} \label{main-thm}
For every translation and scale invariant measure $\mu$ satisfying
$(\star)$, there exists $\lambda_c=\lambda_c(\mu)$, with $0<\lambda_c<\infty$,
such that, with probability one, the complement of the scale invariant soup
with density $\lambda\mu$ contains connected components larger than one point
if $\lambda \leq \lambda_c$, and is totally disconnected if $\lambda>\lambda_c$.
The result holds for the full space soup and for the soup in any domain $D$
with the same $\lambda_c$.
\end{theorem}

We say that a random fractal {\em percolates} if it contains connected
components larger than one point. (This is not the definition
typically used for Mandelbrot's fractal percolation, which involves a
certain crossing event --- see Section~\ref{sec-mandelbrot} below ---
but we think it is more natural and ``canonical,'' at least in the present
context, precisely because it does not involve an arbitrary crossing event.)
Theorem \ref{main-thm} therefore says that for the class of models included
in the statement,
with probability one the system {\em percolates at criticality}.
We will show that this
is equivalent to having positive probability for certain crossing events
involving ``shell-like'' (deterministic) domains. \\

\noindent{\bf Remark} \, As pointed out to us by an anonymous referee,
a standard example of  percolation at criticality is the appearance of a
$k$-ary subtree inside a Galton-Watson tree (e.g., in Bernoulli percolation
on a b-ary tree), with $k \geq 2$. This example has in fact played a role in
fractal percolation (e.g., it is used in the proof of Theorem 1 of \cite{CCD}).
It would be interesting to determine whether there is a connection between
that example and the class of models treated in this paper.

\begin{corollary} \label{cor-large}
Consider a full space soup in ${\mathbb R}^d$ with density $\lambda\mu$,
where $\mu$ is a translation and scale invariant measure satisfying
$(\star)$. If $\lambda \leq \lambda_c(\mu)$, with probability one, the
complement of the soup contains arbitrarily large connected components.
Moreover, if $\mu$ is invariant under rotations, any two open subsets of
${\mathbb R}^d$ are intersected by the same connected component of the
complement of the soup with positive probability.
\end{corollary}

Corollary \ref{cor-large} leaves open the question of existence, and
possibly uniqueness, of an unbounded connected component. We are able
to address this question only for $d=2$ (see Theorem~\ref{thm-uniqueness}). \\

\noindent{\bf Remark} \,
The measure $\mu$ does not need to be completely scale invariant
for our results above to hold. As it will be clear from the proofs,
it suffices that there is an infinite sequence of scale factors
$s_j \downarrow 0$ such that $\mu$ is invariant under scaling by
$s_j$, in the sense described above. This is the case, for instance,
for the multiscale Boolean model studied in Chapter~8
of~\cite{mr} and in~\cite{mpv1,mpv3}.

Indeed, our definition of self-similar
soups is aimed at identifying a natural class of models that is
easy to define and contains interesting examples; we did not try
to define the most general class of models to which our methods
apply.
In Section~\ref{sec-mandelbrot} we will use Mandelbrot's fractal
percolation to illustrate how our techniques can be easily applied
to an even larger class of models.
\vspace{3mm}

The main technical tool in proving Theorem \ref{main-thm} and
Corollary~\ref{cor-large} is Lemma~\ref{main-lemma}, presented
in Section~\ref{sec-proofs}. The lemma implies that the probability
that the complement in a ``shell-like" domain $A$ of a full space
soup contains a connected component that touches both the ``inner"
and the ``outer" boundary of the domain has a discontinuity at some
$0<\lambda_c^A<\infty$, jumping from a positive value at $\lambda_c^A$
to zero for $\lambda>\lambda_c^A$. It is then easy to see that the
complement of the soup must be totally disconnected for
$\lambda>\lambda_c^A$ (Lemma~\ref{lemma2}), which implies that
$\lambda_c^A$ is the same for all ``shell-like" domains and coincides
with the $\lambda_c$ of Theorem~\ref{main-thm}.

For future reference we define what we mean by a {\em shell} and a {\em simple
shell}.
We call a set $A$ a \emph{shell} if it can be written as $A = D \setminus \overline{D'}$,
where $D$ and $D'$ are two non-empty, bounded, $d$-dimensional open sets with
$\overline{D'} \subset D$. A shell $A$ is \emph{simple} if $D$ and $D'$ are
open, concentric ($d$-dimensional) cubes. We will denote by $\Phi_A$ the probability
that the complement of a full space soup contains a connected component that touches
both the ``inner" and the ``outer" boundary of $A$.

We note that the proof of Lemma~\ref{main-lemma} makes essential use of the shell geometry
and would not work in the case, for instance, of crossings of cubes. Throughout the proof
of Theorem~\ref{main-thm}, we choose to use simple shells because they are easier to work
with. However, all our results can be readily generalized to any shell. In particular, we
have the following discontinuity result, which is interesting in its own right.

\begin{corollary} \label{cor}
For all $d \geq 2$, all shells $A$, and all translation and scale
invariant measures $\mu$ satisfying $(\star)$, the following holds:
\begin{itemize}
\item $\Phi_A(\lambda)>0$ if $\lambda \leq \lambda_c(\mu)$,
\item $\Phi_A(\lambda)=0$ if $\lambda > \lambda_c(\mu)$.
\end{itemize}
\end{corollary}

\section{Two-Dimensional Soups} \label{sec-2d}

In two dimensions one can obtain additional information
and show that, like in Mandelbrot's fractal percolation,
a unique infinite connected component appears as soon as
there is positive probability of having connected components
larger than one point (that is, {\it at} and below the
critical point $\lambda_c$).

To prepare for our main result of this section,
Theorem~\ref{thm-uniqueness} below, consider a self-similar
soup in the unit square $(0,1)^2$, and let $g(\lambda)$ be
the probability that the complement of the soup contains a
connected component that crosses the square in the first
coordinate direction, connecting the two opposite sides of
the square. We then have the following result.

\begin{theorem} \label{thm-2D}
For every translation and scale invariant $\mu$ in two dimensions
which satisfies condition $(\star)$ and is invariant under reflections
through the coordinate axes and rotations by 90 degrees, $g(\lambda_c(\mu))>0$.
\end{theorem}

The invariance under reflections through the coordinate axes
is required because $g(\lambda)$ is defined using crossings
of the unit square $(0,1)^2$.

In the case of the Brownian loop soup, Werner states a version of
Theorem~\ref{thm-2D} in~\cite{werner1,werner2}. The choice of the
unit square in Theorem~\ref{thm-2D} is made only for convenience,
and similar results can be proved in the same way for more general
domains. The reflection invariance is a technical condition needed
in the proof in order to apply a technique from~\cite{dm} (see the
proof of Lemma~5.1 there). The same technique, combined with
Theorem~\ref{thm-2D}, can be used to prove the next theorem, which
is our main result of this section.

\begin{theorem} \label{thm-uniqueness}
For every translation and scale invariant $\mu$ in two dimensions which
satisfies condition $(\star)$ and is invariant under reflections through
the coordinate axes and rotations by 90 degrees, if $\lambda \leq \lambda_c(\mu)$,
the complement of the full plane soup with density $\lambda\mu$ has a
unique unbounded component with probability one.
\end{theorem}

The informed reader might believe that the uniqueness result would follow
from a version of the classical Burton-Keane argument (see \cite{BK}).
However, in the Burton-Keane argument it is crucial, for instance, that a path
from the inside to the outside of a cube of side length $n$ uses at least
(roughly) a portion $1/{n}^{d-1}$ of the ``surface'' of the cube (e.g.,
the number of sites of ${\mathbb Z}^d$ on the boundary,  for a lattice
model defined on ${\mathbb Z}^d$), so that there is enough space for at most
$O(n^{d-1})$ disjoint paths. There is clearly no analogue of this for the
continuous models in this paper (nor for Mandelbrot percolation), since the
relevant paths have no ``thickness.''

\section{Applications to Mandelbrot's Fractal Percolation} \label{sec-mandelbrot}

The method of proof of Theorem~\ref{main-thm} works in greater generality
than the class of scale invariant soup models introduced in this paper. In
order for the method to work, it suffices to have some form of scale invariance.
To illustrate this fact, we will consider a well-known model, called
\emph{fractal percolation}, that
was introduced by Mandelbrot~\cite{mandelbrot,mandelbrot-book} and is defined by
the following iterative procedure.

For any integers $d \geq 2$ and $N \geq 2$, and real number $0<p<1$,
one starts by partitioning the unit cube $[0,1]^d \subset {\mathbb R}^d$
into $N^d$ subcubes of equal size. Each subcube is independently retained
with probability $p$ and discarded otherwise. This produces a random set
${\mathcal C}_N^1={\mathcal C}_N^1(d,p) \subset [0,1]^d$. The same procedure
is then repeated inside each retained subcube, generating the random set
${\mathcal C}_N^2 \subset {\mathcal C}_N^1$. Iterating the procedure ad
infinitum yields an infinite sequence of random sets
$[0,1]^d \supset \ldots \supset {\mathcal C}_N^k \supset {\mathcal C}_N^{k+1} \supset \ldots$ .
It is easy to see that the \emph{limiting retained set}
${\mathcal C}_N:=\cap_{k=1}^{\infty} {\mathcal C}_N^k$ is well defined.

Several authors studied various aspects of Mandelbrot's fractal percolation,
including the Hausdorff dimension of ${\mathcal C}_N$, as detailed in~\cite{DG},
and the possible existence of paths
\cite{CCD,dm,Meester,CCGS,white,ChCh,FG1,FG2,bc} and
$(d-1)$-dimensional ``sheets" \cite{CCGS,Orzechowski,bc} traversing the unit cube
between opposite faces. Dekking and Meester~\cite{dm} proposed a ``morphology"
of more general ``random Cantor sets,'' obtained by generalizing the successive
``deletion of middle thirds'' construction using random substitutions,
and showed that there can be several critical points at which the connectivity
properties of a set change. Accounts of fractal percolation can be found
in~\cite{lchayes} and Chapter~15 of~\cite{falconer-book}.

In this section we define three potentially different critical points.
Theorem~\ref{thm-mandelbrot} shows that two of them are in fact the same.
Furthermore, we prove that the third one is equal to the other two for $N$
large enough, and conjecture that they are in fact the same for all $N$.

The first critical point is
\[
\tilde{p}_c = \tilde{p}_c(N,d) :=
\sup \{ p: {\mathcal C}_N \text{ is totally disconnected with probability one} \}.
\]
To define the second critical point we focus on a specific shell.
This choice is convenient but arbitrary and unnecessarily restrictive.
Indeed, the proof of the next theorem shows that we could have chosen any
other shell, so that $\hat{p}_c$ defined below is independent of the choice
of shell. Let $A \subset [0,1]^d$ be the domain obtained by removing from
the open unit cube $(0,1)^d$ the cube $(1/2,\ldots,1/2)+[0,1/3]^d$ of side
length $1/3$, centered at $(1/2,\ldots,1/2)$. Denote by $\phi_A(p)$ the
probability that the limiting retained set contains a connected component
that intersects both the ``inner'' and the ``outer'' boundary of $A$. The second
critical point is $\hat{p}_c=\hat{p}_c(N,d) := \inf \{p : \phi_A(p)>0 \}$.
Our first result of this section concerns $\hat{p}_c(N,d)$ and $\tilde{p}_c(N,d)$.

\begin{theorem} \label{thm-mandelbrot}
For all $d \geq 2$ and $N \geq 2$, $\hat{p}_c=\hat{p}_c(N,d)$
satisfies $0<\hat{p}_c<1$. Moreover $\phi_A(\hat{p}_c)>0$, while
${\mathcal C}_N$ is a.s.\ totally disconnected when $p<\hat{p}_c$.
Hence, $\hat{p}_c(N,d)=\tilde{p}_c(N,d)$ for every $N$ and $d$.
\end{theorem}

Mandelbrot's fractal percolation can be extended to a full space model
by tiling ${\mathbb R}^d$ with independent copies of the system in the natural
way. We call this model \emph{full space fractal percolation}. As a consequence
of the previous theorem we have the following result.

\begin{corollary} \label{cor-mandelbrot}
Consider full space fractal percolation with $d \geq 2$ and $N \geq 2$. With
probability one, the limiting retained set contains arbitrarily large connected
components for $p \geq \hat{p}_c$, and is totally disconnected for $p<\hat{p}_c$.
\end{corollary}

We say that there is a (left to right) crossing of the unit cube if ${\mathcal C}_N$
contains a connected component that intersects both $\{0\}\times[0,1]^{d-1}$ and
$\{1\}\times[0,1]^{d-1}$. Let $p_c(N,d)$ be the infimum over all $p$ such that
there is a crossing of the unit cube with positive probability. Sometimes the
system is said to percolate when such a crossing occurs.
For $d=2$ and all $N \geq 2$, Chayes, Chayes and Durrett~\cite{CCD} discovered
that, at the critical point $p_c(N,2)$, the probability of a crossing is strictly
positive (see~\cite{dm} for a simple proof). A slightly weaker result in three
dimensions was obtained in~\cite{CCGS}. Broman and Camia~\cite{bc} were able to
extend the result of Chayes, Chayes and Durrett to all $d \geq 2$, but only for
sufficiently large $N$. However, the same is conjectured to hold for all $N$.

It is interesting to notice that in two dimensions one can prove that, for $p=p_c(N,2)$,
${\cal C}_N$ contains an infinite connected component with probability one~\cite{CCD}.
This is in sharp contrast with lattice percolation, where it has been proved that,
with probability one, the system does not have an infinite cluster at the critical
point in dimensions 2 and $\geq 19$. (The same is believed to hold in all dimensions
--- see~\cite{grimmett-book} for a general account of percolation theory.)

By Theorem~\ref{thm-mandelbrot},
${\mathcal C}_N$ is totally disconnected with probability one when $p<\hat{p}_c$,
so that $\hat{p}_c(N,d) \leq p_c(N,d)$ for all $N$ and $d$.
Furthermore we have the following result.
\begin{theorem} \label{thm-equivalence}
For $d=2$ we have that $\hat{p}_c(N,2)=p_c(N,2)$ for all $N \geq 2$.
Furthermore, for every $d \geq 3$, there exists $N_0=N_0(d)$ such that, for all $N \geq N_0$,
$\hat{p}_c(N,d)=p_c(N,d)$.
\end{theorem}

\noindent{\bf Remark} \, We conjecture that $\hat{p}_c(N,d)=p_c(N,d)$ for all $N \geq 2$
and all $d \geq 2$.

\section{Proofs} \label{sec-proofs}

\subsection{Proofs of the Main Results}

Before we can give the actual proofs, we need some definitions.
Let $A = D \setminus \overline{D'}$ be a shell.
Given a set $K$, if $A \setminus K$ contains a connected
component that connects the boundary of $D$ with that of $D'$ (in
other words, if $K$ does not disconnect $\partial D'$ from $\partial D$),
we say that the complement of $K$ \emph{crosses} $A$, or that there is a
crossing of $A$ in $A \setminus K$. We let $\Phi_A(\lambda)$ denote the
probability that the complement of a full space soup crosses $A$.
If $\mu$ satisfies condition $(\star)$, the function $\Phi_A(\lambda)$
is left-continuous in $\lambda$.
To see this, consider $\Phi_A^{\varepsilon}(\lambda)$,
the analogous crossing probability obtained by disregarding sets in the soup
of diameter smaller than $\varepsilon$. A standard coupling of Poisson processes
(between different values of $\lambda$) shows that, if condition $(\star)$ is
satisfied, $\Phi_A^{\varepsilon}(\lambda)$ is continuous in $\lambda$ for any
$\varepsilon>0$. Furthermore, $\Phi_A(\lambda)$ corresponds to the $\varepsilon \to 0$
limit of $\Phi_A^{\varepsilon}(\lambda)$, and is therefore left-continuous in $\lambda$,
since $\Phi_A^{\varepsilon}(\lambda)$ is nonincreasing in $\lambda$.

We now start with the proof of our main theorem, leaving out some lemmas
that will be proved later.

\bigskip

\noindent{\bf Proof of Theorem~\ref{main-thm}}

\bigskip

\noindent{\bf Full Space Soup} \,
We will first prove the result for the full space soup. Define $\lambda_c=\lambda_c(\mu)$ to
be the infimum of all $\lambda$ such that with probability one
the complement of the full space soup contains at most isolated points.
Clearly, $\Phi_A(\lambda)=0$ for
$\lambda>\lambda_c$. The following lemma, whose proof is standard
and deferred till later on, holds.

\begin{lemma} \label{lemma1}
For any translation and scale invariant measure $\mu$ satisfying
condition $(\star)$, we have $0<\lambda_c(\mu)<\infty$.
\end{lemma}

In order to conclude the proof for the full space soup, it suffices to
show that $\Phi_A(\lambda_c)>0$ for some simple shell $A$,
since that would imply that the complement of the soup cannot be totally
disconnected with probability one at $\lambda=\lambda_c$. In order to
achieve that, we combine the left-continuity of $\Phi_A(\lambda)$ with the
two following lemmas.

\begin{lemma} \label{lemma2}
If $\lambda$ is such that $\Phi_A(\lambda)=0$ for some simple shell $A$,
then the complement of the full space soup with density $\lambda\mu$
is totally disconnected with probability one.
\end{lemma}

\begin{lemma} \label{main-lemma}
For any simple shell $A$, there exists an $\varepsilon>0$
such that, if $\Phi_A(\lambda) \leq \varepsilon$, then $\Phi_A(\lambda)=0$.
\end{lemma}

Lemma~\ref{lemma2} implies that for all simple shells $A$,
$\Phi_A(\lambda)>0$ for $\lambda<\lambda_c$. Together with with the
left-continuity of $\Phi_A(\lambda)$ and Lemma~\ref{main-lemma}, this
implies that $\Phi_A(\lambda_c)>0$, which concludes this part of the proof.

\bigskip

\noindent{\bf Soup in a Bounded Domain} \, We now prove the result
for the soup in a domain $D$. Let $\lambda_c^D$ be the infimum
of the set of $\lambda$'s such that the complement of the soup with
intensity $\lambda\mu$ in $D$ is totally disconnected with probability
one. Coupling the soup in $D$ with a full space soup with cutoff larger
than the maximum radius allowed in $D$ by using the same Poisson
realization for both before applying the cutoff or the condition that
sets be contained in $D$, one can easily see that $\lambda_c^D\geq\lambda_c$.
Indeed, more sets are ``discarded" in the case of the soup in $D$, meaning
that the intersection with $D$ of the complement of the full space soup
is contained in the complement of the soup in $D$. For any $\lambda<\lambda_c$,
the complement of the full space soup intersected with $D$ contains a connected
component larger than one point with positive probability. This is because we
can cover ${\mathbb R}^d$ with a countable number of copies of $D$ and use
translation invariance. It follows that the complement of the soup in $D$,
must also contain a connected component larger than one point with positive
probability showing that $\lambda\leq\lambda_c^D$ so that $\lambda_c^D\geq\lambda_c$.

On the other hand, for any closed set $G \subset D$, when $\lambda>\lambda_c$,
the intersection with $G$ of the complement of the soup in $D$ is easily seen
to be totally disconnected by comparing it with the intersection of $G$ with
the complement of a full space soup with cutoff smaller than
$\frac{1}{2}\text{dist}(G,\partial D)$, coupled to the soup in $D$ in the
same way as before. Therefore $\lambda_c^D\leq\lambda_c$, and we conclude
that $\lambda_c^D=\lambda_c$ for all domains $D$.

It remains to check that at the critical point $\lambda_c$, the complement
of the soup in $D$ contains connected components larger than one point.
This can be done by coupling the soup in $D$, as before, to a full space
soup with cutoff larger than the maximum radius allowed in $D$. The
intersection with $D$ of the complement of such a soup is contained in the
complement of the soup in $D$. For $\lambda=\lambda_c$, the complement of
the full space soup intersected with $D$ contains a connected component
larger than one point with positive probability, and so does the complement
of the soup in $D$. The proof of Theorem~\ref{main-thm} is therefore complete. \fbox{} \\

We now turn to the proofs of Lemmas \ref{lemma1}, \ref{lemma2} and \ref{main-lemma}.
Lemma~\ref{lemma1} can be proved in various, rather standard, ways. A detailed
proof in the context of the multiscale Boolean model can be found in
Chapter~8 of~\cite{mr}, while a different proof in the context of the Brownian
loop soup is sketched in~\cite{werner2}. Both proofs are given in two dimensions,
but the dimensionality of the space is irrelevant and the same arguments work
in all dimensions. Since the same ideas work for all scale invariant soups, we only
sketch the proof of the lemma, and refer the interested reader to~\cite{mr}
or~\cite{werner2} for more details.

\bigskip

\noindent{\bf Proof of Lemma~\ref{lemma1}} \,
Since the whole space can be partitioned in a countable number of cubes, in
order to prove that $\lambda_c<\infty$, it suffices to show that, for some
$\lambda$ sufficiently large, the unit cube $[0,1]^d$ is completely covered
with probability one.

Let $n$ be a positive integer. Given $\mu$, choose $\alpha$ small enough so
that, with strictly positive probability, the cube $(-\alpha, \alpha)^d$
is covered by a set from the soup of intensity $\mu$ contained inside the cube
$(-1, 1)^d$.
Note that, if the probability of the event described above is strictly positive
for a soup of intensity $\mu$, and if we denote by $a(\lambda)$ the probability
of the same event for a soup of intensity $\lambda\mu$, we have $a(\lambda) \to 1$
as $\lambda \to \infty$.

We can cover the unit cube with order $2^n$ cubes of side length $\alpha 2^{-n+1}$.
Each such cube is contained inside order $n$ nested simple shells with diameter
$\leq 2\sqrt{2}$ and constant ratio $\alpha$ between the side length of the
outer cube and that of the inner cube. Therefore,  by scale invariance the probability
that the unit cube is not covered by sets from a soup with intensity $\lambda\mu$ is bounded
above by some constant times
\begin{equation} \nonumber
 2^n (1 - a(\lambda))^n = 2^{(1 - |\log_2(1 - a(\lambda))|)n}.
\end{equation}
For $\lambda$ so large that $|\log_2(1 - a(\lambda))| > 1$,
the exponent in the bound is negative and so the bound tends to zero as
$n \to \infty$, which concludes the proof of the first part of the lemma.

Now let $\cal K$ denote the collection of sets (from a soup) that intersect the
unit cube $[0,1]^d$. For each set $K \in {\cal K}$, define $q(K) \in {\mathbb N}$
in such a way that $\text{diam}(K) \in (2^{-q(K)-1}, 2^{-q(K)}]$. Partition the unit cube
in cubes of side length $2^{-n}$, and for each cube $C$ of the partition $\Pi_n$ define
\begin{equation} \nonumber
\tilde X(C) = 1_{\{\not \exists K : q(K)=n \text{ and } K \cap C \neq \emptyset\}}.
\end{equation}
The random variables $\{ \tilde X(C) \}_{C \in \Pi_n}$ are not independent, since
two adjacent elements of $\Pi_N$ can intersect the same set $K$. It is in fact
easy to see that they are positively correlated. However, 
it is possible to couple the collection of random variables $\{ \tilde X(C) \}_{C \in \Pi_n}$
with a collection of independent Bernoulli random variables $\{ X(C) \}_{C \in \Pi_n}$ such that,
for each $C$, $\tilde X(C) \geq X(C)$ almost surely (see, e.g., \cite{lss}). Moreover,
for any $\delta>0$ we can take $P(X(C)=1) \geq 1- \delta$, if $P(\tilde X(C)=1)$
is close enough to 1.

Let $F_n(C)$ denote the collection of compact sets with nonempty interior that
intersect $C$ and have diameters in $(2^{-n-1}, 2^{-n}]$. Then, for all $C \in \Pi_n$,
\begin{equation} \nonumber
P(\tilde X(C)=1) = e^{-b\lambda},
\end{equation}
where $b=\mu(F_n(C))<\infty$ by scale (and translation) invariance and condition $(\star)$.
Therefore, by taking $\lambda$ small enough, we can make $P(\tilde X(C)=1)$, and thus
also $P(X(C)=1)$, arbitrarily close to 1. Since the collection of random variables
$\{ X(C) : C \in \Pi_n, n \in {\mathbb N} \}$ defines a Mandelbrot percolation process
in $[0,1]^d$ with retention probability equal to $P(X(C)=1)$, and whose retained set
is contained in $[0,1]^d \setminus {\cal K}$, for sufficiently small $\lambda$,
$[0,1]^d \setminus {\cal K}$ contains connected components larger than one
point with positive probability. \fbox{} \\

\noindent{\bf Proof of Lemma~\ref{lemma2}} \,
Let $\lambda>0$ and the simple shell $A=B^{out} \setminus \overline{B^{in}}$
be such that $\Phi_A(\lambda)=0$. Because of scale and translation
invariance, $\Phi_{A'}(\lambda)=0$ for any $A'$ obtained by translating
a scaled shell $sA$ with $s \leq 1$.

Given $\varepsilon>0$, take $s=s(\varepsilon)$ such that $0<s<1$ and
$\text{diam}(sB^{out})=s \, \text{diam}(B^{out})<\varepsilon$, and consider the
simple shell $sA$. Consider a tiling of ${\mathbb R}^d$ with non-overlapping
(except along the boundaries) translates of $\overline{sB^{in}}$ such that the
centers of the cubes form a regular lattice isomorphic to ${\mathbb Z}^d$. If
the full space soup contains a connected component of diameter larger than
$\varepsilon$, such a component must intersect some of the cubes from the tiling.

For any cube from the tiling, the probability that it is intersected by
a connected component of the complement of the soup of diameter larger than
$\varepsilon$ is bounded above by the probability that a translate of $sA$
is crossed by the complement of the soup. This follows from the fact that
the diameter of the connected component is strictly larger than the diameter
of $sA$. Since the probability of crossing $sA$ is
zero, and the number of cubes in the tiling is countable, we conclude
that the full space soup cannot contain a connected component of diameter
larger than $\varepsilon$, for any $\varepsilon>0$. \fbox{} \\

\noindent{\bf Proof of Lemma~\ref{main-lemma}} \,
Let $B^{in}$ and $B^{out}$ be the two $d$-dimensional, concentric, open cubes
such that $A=B^{out} \setminus \overline{B^{in}}$.
For $0<s<1$, consider a tiling of ${\mathbb R}^d$
with non-overlapping (except along the boundaries) translates of $sB^{in}$ such that the centers of the cubes form
a regular lattice isomorphic to ${\mathbb Z}^d$. We can use this isomorphism
to put the translates of $sB^{in}$ in a one-to-one correspondence with the
vertices of ${\mathbb Z}^d$, and thus index them via the vertices of ${\mathbb Z}^d$.

For each translate $B^{in}_{x,s}$, $x \in {\mathbb Z}^d$, of $sB^{in}$, consider
the translate $B^{out}_{x,s}$ of $sB^{out}$ concentric to $B^{in}_{x,s}$. The two
define the simple shell
$A_{x,s} = B^{out}_{x,s} \setminus \overline{B^{in}_{x,s}}$.
We then have a collection $\{ A_{x,s} \}_{x \in {\mathbb Z}^d}$ of
(overlapping) simple shells indexed by ${\mathbb Z}^d$.

Let ${\cal K}_s$ denote the collection of sets from the full space soup with
diameter at most $s$. Obviously, ${\cal K}_s$ is distributed like
a full space soup with cutoff $\delta=s$. Let $K_s=\bigcup_{K \in {\cal K}_s}K$,
and denote by $\Psi_x(\lambda,s)$ the probability that there is a crossing of
$A_{x,s}$ in the complement of $K_s$. It immediately follows from scale and
translation invariance and the way $A_{x,s}$ has been defined that
$\Psi_x(\lambda,s)=\Phi_{s^{-1}A_{x,s}}(\lambda) = \Phi_A(\lambda)$.

We now introduce the graph ${\mathbb M}^d$ whose set of vertices is ${\mathbb Z}^d$
and whose set of edges is given by the adjacency relation: $x \sim y$ if and only
if $||x-y||=1$, where $||x||=||(x_1,\ldots,x_d)||:=\max_{1 \leq i \leq d} |x_i|$
and $|\cdot|$ denotes absolute value.
Next, for each $0<s<1$, we define the random variables $\{ X_s(x) \}_{x \in {\mathbb Z}^d}$
by letting $X_s(x)=1$ if there is a crossing of $A_{x,s}$ in the complement of $K_s$,
and $X_s(x)=0$ otherwise. By construction, the probability that $X_s(x)=1$ equals
$\Psi_x(\lambda,s)=\Phi_{s^{-1}A_{x,s}}(\lambda) = \Phi_A(\lambda)<1$.

Note that, if $||x-y||>[\text{diam}(B^{out})+2]/l$, where $l$ is the Euclidean
side length of $B^{in}$, then $X_s(x)$ and $X_s(y)$ are independent of each other.
This implies that we can apply Theorem~B26 of~\cite{SIS} (see p.~14 there; the
result first appeared in~\cite{lss}) to conclude that there exist i.i.d.\ random
variables $\{ Y_s(x) \}_{x \in {\mathbb Z}^d}$ such that $Y_s(x)=1$ with probability
$p<1$ and $Y_s(x)=0$ otherwise, and $Y_s(x) \geq X_s(x)$ for every $x\in{\mathbb Z}^d$.
Moreover, one can let $p \to 0$ as $\Phi_A(\lambda) \to 0$.

For each $0<s<1$, using the random variables $\{ Y_s(x) \}_{x \in {\mathbb Z}^d}$,
we can define a Bernoulli site percolation model on ${\mathbb M}^d$ by declaring
$x \in {\mathbb Z}^d$ open if $Y_s(x)=1$ and closed if $Y_s(x)=0$. We denote
by $p_c(d)$ the critical value for Bernoulli site percolation on ${\mathbb M}^d$.
(See~\cite{grimmett-book} for a general account on percolation theory.)

Let $G_s := \{x \in {\mathbb Z}^d : A_{x,s} \subset A \}$, i.e., the set of vertices
of ${\mathbb Z}^d$ corresponding to simple shells $A_{x,s}$ contained in $A$. Note
that, if $s$ is sufficiently small, ${\mathbb Z}^d \setminus G_s$ contains two components,
of which one is unbounded. These components are connected in terms of the adjacency relation
$\sim$ used to define ${\mathbb M}^d$ when considered as subsets of the vertex set of ${\mathbb M}^d$.

When this is the case, if there is a crossing of $A$ in the complement of the full space
soup, then the percolation process on ${\mathbb M}^d$ defined above has an open cluster
that connects the bounded component to the unbounded component of ${\mathbb Z}^d \setminus G_s$,
``crossing'' $G_s$. The reason is that if the crossing of $A$ intersects a box $B^{in}_{x,s}$,
then $A_{x,s}$ must be also be crossed and so $Y_s(x)=1$.
The diameter of such an open cluster is at least of order
$\text{dist}(\partial B^{out}, \partial B^{in})/sl$, and the cluster is contained in $G_s$, whose
diameter is of the order of $\text{diam}(B^{out})/sl$. (Here, for $K_1,K_2 \subset {\mathbb R}^d$,
$\dist(K_1,K_2):=\inf\{|x-y| : x \in K_1, y \in K_2 \}$.)

We are now ready to conclude the proof. Assume that $\Phi_A(\lambda)<\varepsilon$ with
$0<\varepsilon<1$ so small that one can choose $p=P(Y_s(x)=1)$ so that $p<p_c(d)$. Take
$s$ so small that $G_s$ contains two connected components, as explained above. Then,
for every $s$ sufficiently small, $\Phi_A(\lambda)$ is bounded above by the probability
that the Bernoulli percolation process defined via the random variables
$\{ Y_s(x) \}_{x \in {\mathbb Z}^d}$ contains an open cluster of diameter at least $L/s$
inside a region of linear size at most $L'/s$, for some $L,L'<\infty$. Since $p<p_c(d)$,
it follows from standard percolation results (see, e.g., \cite{grimmett-book}) that the
probability of such an event goes to zero as $s \to 0$, proving the lemma. \fbox{} \\

\noindent{\bf Proof of Corollary~\ref{cor}} \,
Let $A$ be a shell and let $A'$ be a simple shell such that $A\subset A'$.
Using Theorem \ref{main-thm} we conclude that $\Phi_A(\lambda)=0$ if $\lambda > \lambda_c(\mu)$.
Furthermore, we have that
\[
0<\Phi_{A'}(\lambda_c(\mu))\leq \Phi_A(\lambda_c(\mu)),
\]
since a crossing of $A'$ implies a crossing of $A$. \fbox{} \\

\bigskip

\noindent{\bf Remark} \, It is possible to define the notion of $(d-1)$-dimensional
crossings of general shells. For example, a crossing could be, informally, a connected
subset of the complement of the soup which divides the shell into two disjoint
parts, both touching the ``inner'' and ``outer'' boundary of the shell. Using
this definition of crossing it is possible to show results analogous to Theorem \ref{main-thm}
and Corollary \ref{cor}. Note that, for $d=2$ this is not the definition that we use, but
it is easy to see that our results would still be true with this definition of crossing.
\vspace{3mm}

The proof of our second main result is now easy.

\bigskip

\noindent{\bf Proof of Corollary~\ref{cor-large}} \, The first claim can
be proved using Theorem~\ref{main-thm} and a simple scaling argument, but it
is also an immediate consequence of Corollary~\ref{cor} combined with
translation invariance.

The proof of the second claim uses rotation invariance. Let us consider,
without loss of generality, two disjoint open balls, $B_1$ and $B_2$,
of radii $r_1$ and $r_2$ and centered at $x_1$ and $x_2$, respectively.
(If the balls are not disjoint, the complement of the soup intersects
$B_1 \cap B_2$ with a connected component larger than one point with
positive probability.) Let $d_{12}=|x_1-x_2|$, and consider the
shell $A = \{ x : |x-x_1|<d_{12} \} \setminus \{ x : |x-x_1| \leq r_1/2 \}$.
From Corollary~\ref{cor} we know that $\Phi_A(\lambda)>0$ for $\lambda \leq \lambda_c$.
It then follows from rotation invariance that there is positive probability
that the complement of the soup contains a connected component that connects
the sphere $\{ x : |x-x_1| = r_1/2 \}$ with the surface
$\{ x : |x-x_1| = d_{12} \} \cap B_2$. Such a connected component must
intersect both $B_1$ and $B_2$. \fbox{} \\

\bigskip

\subsection{Proofs of the Additional Two-Dimensional Results}

\noindent{\bf Proof of Theorem~\ref{thm-2D}} \, Let $\tilde{f}(\lambda)$
denote the probability that the complement of the full plane soup with density
$\lambda\mu$ and cutoff $\delta=2$ contains a connected component that crosses
the rectangle $[0,1] \times [0,2]$ horizontally, and $\tilde{g}(\lambda)$ the
probability that it contains a connected component that crosses the square
$[0,1]^2$ horizontally. Consider the annulus $A = [-3/2,3/2]^2 \setminus [-1/2,1/2]^2$.
By Corollary \ref{cor}, $\Phi_A(\lambda)>0$ if $\lambda \leq \lambda_c(\mu)$.
It is easy to see (Figure~\ref{fig4}) that any crossing of $A$ must cross either
a square of side length $1$ or a rectangle of side lengths $1$ and $2$ in the
``easy'' direction. Using translation and rotation invariance, this implies
$\tilde{f}(\lambda_c)>0$.

Let us now couple the full plane soup with cutoff $\delta=2$ with full plane fractal
percolation with $N=3$ in such a way that at level $k=0,1,\ldots$ of the fractal
percolation construction, a square of side length $3^{-k}$ is discarded if and only
if it is covered by sets of the soup of diameter between $2/3^k$ and $2/3^{k+1}$,
which happens with positive probability $q$. It is immediate that the limiting
retained set of the full plane fractal percolation process contains the complement
of the full plane soup. Therefore, $\tilde{f}(\lambda_c)>0$ implies that there is
positive probability that the limiting retained set of the fractal percolation process
contains a connected component that crosses the rectangle $[0,1] \times [0,2]$ horizontally.

Note that in the fractal percolation process defined above, squares are not discarded
independently, due to the presence of sets that can intersect two or more squares
(up to four). However, two level-$k$ squares of side length $3^{-k}$ at distance
larger than $2/3^k$ are retained or discarded independently. In particular, if one
marks every third square in a line of level-$k$ squares, all marked squares are
discarded with probability $q>0$, independently of each other. This observation
implies that we can apply the proof of Lemma~5.1 of~\cite{dm} to the fractal
percolation process defined above (as the reader can easily check).

The proof of the lemma shows that horizontal crossings of the rectangle
$[0,1] \times [0,2]$ cannot be ``too straight,'' they must possess a certain
``wavyness'' so that, using invariance under reflections through the $y$-axis
and translations, and the fact that crossing events are positively correlated,
the horizontal crossings in five partially overlapping $1 \times 2$ rectangles
``hook up'' with positive probability to form a horizontal crossing of the
rectangle $[0,3] \times [0,2]$ (see Figure~5 of~\cite{dm} and the discussion
in the proof of Lemma~5.1 there).

Since the horizontal crossings of the rectangle $[0,1] \times [0,2]$ in the complement
of the soup form a subset of the fractal percolation crossings, they must possess the
same ``wavyness'' property. In our setting, positive correlation of crossing events
follows, for instance, from~\cite{janson} and the fact that crossing events for the
complement of the soup are decreasing.
(Let ${\cal K} \subset {\cal K}'$ denote two soup realizations; an event $\cal A$ is
decreasing if ${\cal K} \notin {\cal A}$ implies ${\cal K}' \notin {\cal A}$.)
Therefore, using the same ``hook up'' technique as in the proof of Lemma~5.1 of~\cite{dm},
but with crossings in the complement of the soup, we can conclude that there is positive
probability that the complement of the full plane soup contains a horizontal crossing of
the rectangle $[0,3] \times [0,2]$, and thus also of the square $[0,3] \times [0,3]$. By
scaling, this implies that $\tilde{g}(\lambda_c)>0$.

To conclude the proof, we couple the soup in the unit square $(0,1)^2$ with density
$\lambda_c(\mu)\mu$ with the full space soup with the same density and cutoff $\delta=2$
by using the same Poisson realization for both before applying the cutoff or the condition
that discs be contained in $(0,1)^2$. Clearly, the intersection with $(0,1)^2$ of the
complement of the full space soup is contained in the complement of the soup in $(0,1)^2$.
Therefore, $\tilde{g}(\lambda_c)>0$ implies $g(\lambda_c)>0$, as required. \fbox{} \\

\noindent{\bf Proof of Theorem~\ref{thm-uniqueness}} \,
The proof of Theorem~\ref{thm-2D} shows that, for all $\lambda \leq \lambda_c(\mu)$,
there is positive probability that the complement of the full plane soup with density
$\lambda\mu$ and cutoff $\delta=1$ contains a horizontal crossing of the rectangle
$[0,3] \times [0,2]$. Crossing events like the one just mentioned are decreasing and
are therefore positively correlated (see, e.g., Lemma 2.2 of~\cite{janson}).

\begin{figure}[!ht]
\begin{center}
\includegraphics[width=9cm]{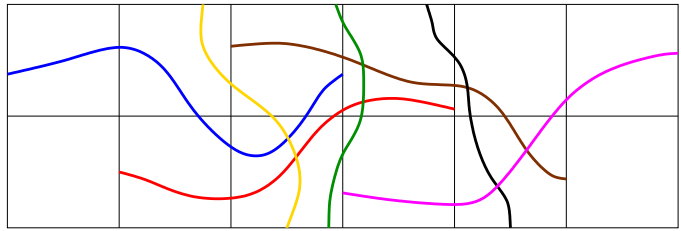}
\caption{Pasting horizontal crossings of 3 by 2 rectangles and vertical crossings
of 2 by 2 squares, a crossing of a 6 by 2 rectangle is obtained. (The different
colors for the crossings serve only to enhance the visibility of the figure.)}
\label{fig1}
\end{center}
\end{figure}

Let us call $a(\lambda)$ the probability that the complement of the full plane soup with
density $\lambda\mu$ and cutoff $\delta=1$ contains a horizontal crossing of the rectangle
$[0,3] \times [0,2]$, and $b(\lambda)$ the probability that it contains a vertical crossing
of the square $[0,2] \times [0,2]$. Positive correlation of crossing events, combined with
a standard ``pasting'' argument (see Fig.~\ref{fig1}), implies that the probability that
the complement of the full plane soup with intensity $\lambda\mu$ and cutoff $\delta=1$
contains a horizontal crossing of the rectangle $[0,6] \times [0,2]$ is bounded below
by $a(\lambda)^4 b(\lambda)^3$. We denote by $h_0$ this probability, and by $h_n$ the
probability of the event ${\cal B}_n$ that the complement of the soup contains a
horizontal crossing of $[0,2 \cdot 3^{n+1}] \times [0,2 \cdot 3^n]$.

Consider the full plane soup with cutoff $\delta_n=3^{-n}$ obtained by a ``thinning''
of the soup with cutoff $\delta=1$ that consists in removing from it all the sets
with diameter larger than $3^{-n}$. By scaling, the probability that the complement
of the soup with cutoff $\delta_n$ contains a horizontal crossing of the rectangle
$[0,3] \times [0,2]$ is equal to $h_n$. The complements of the soups with cutoffs
$\{\delta_n\}_{n \in {\mathbb N}}$ form an increasing (in the sense of inclusion
of sets) sequence of nested sets. Therefore, the limit of $h_n$ as $n \to \infty$
is the probability of $\bigcup_{n \geq 0} {\cal B}_n$. By Kolmogorov's zero-one law,
the latter probability is either 0 or 1. However, since it cannot be smaller than
$h_0>0$, it must necessarily be 1.

\begin{figure}[!ht]
\begin{center}
\includegraphics[width=10cm]{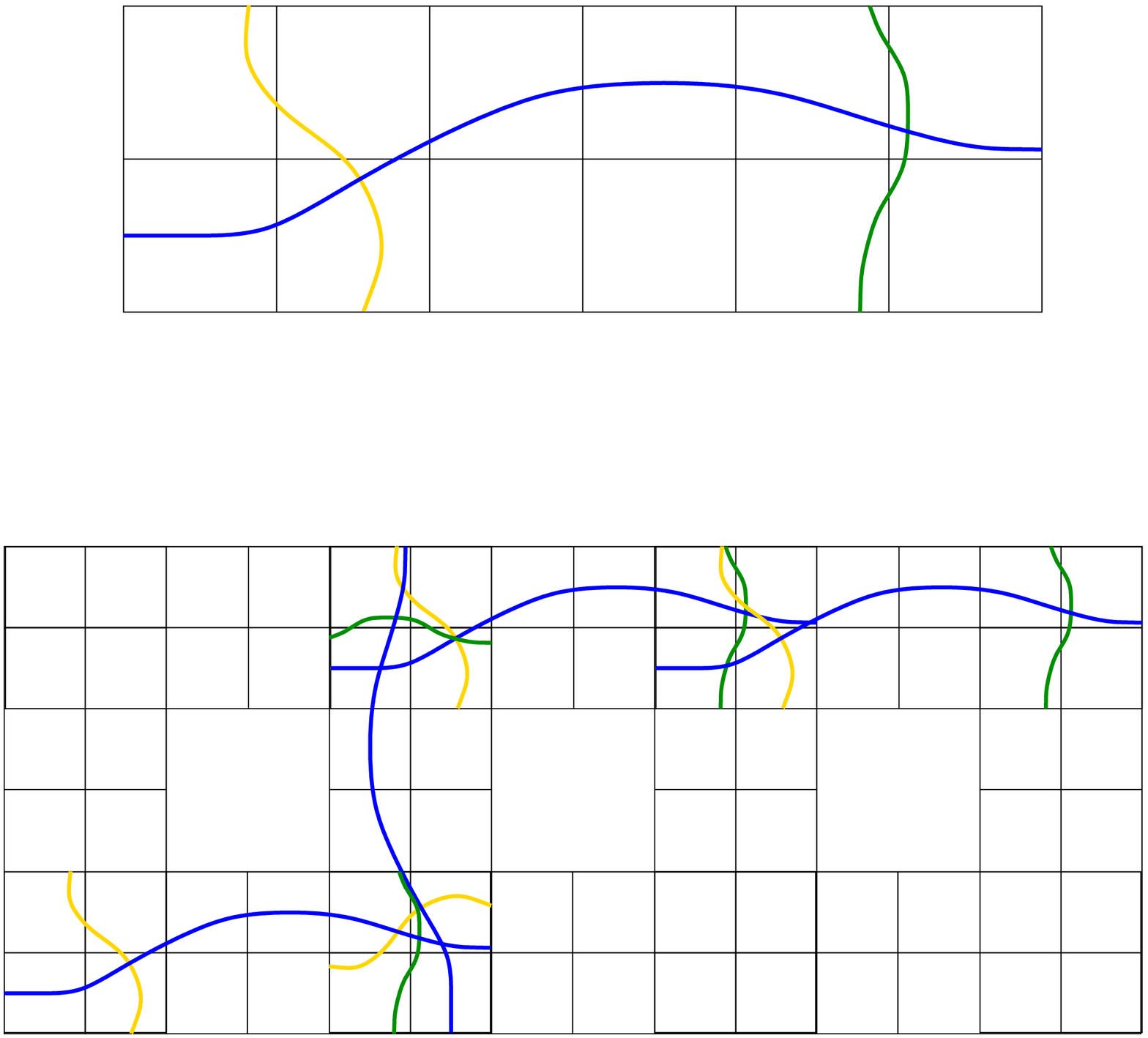}
\caption{The event depicted above and its rotation by 90 degrees can be used to
couple the complement of the soup to a one-dependent percolation process on a
square grid whose open edges correspond to rectangles where the event occurs.
We denote by $p_n$ the probability of the event when the elementary squares in
the figure have side length $3^n$, i.e., the probability that an edge is open
in the corresponding bond percolation process.(As in Fig.~\ref{fig1}, the different
colors for the crossings serve only to enhance the visibility of the figure.)}
\label{fig2}
\end{center}
\end{figure}

\begin{figure}[!ht]
\begin{center}
\includegraphics[width=5cm]{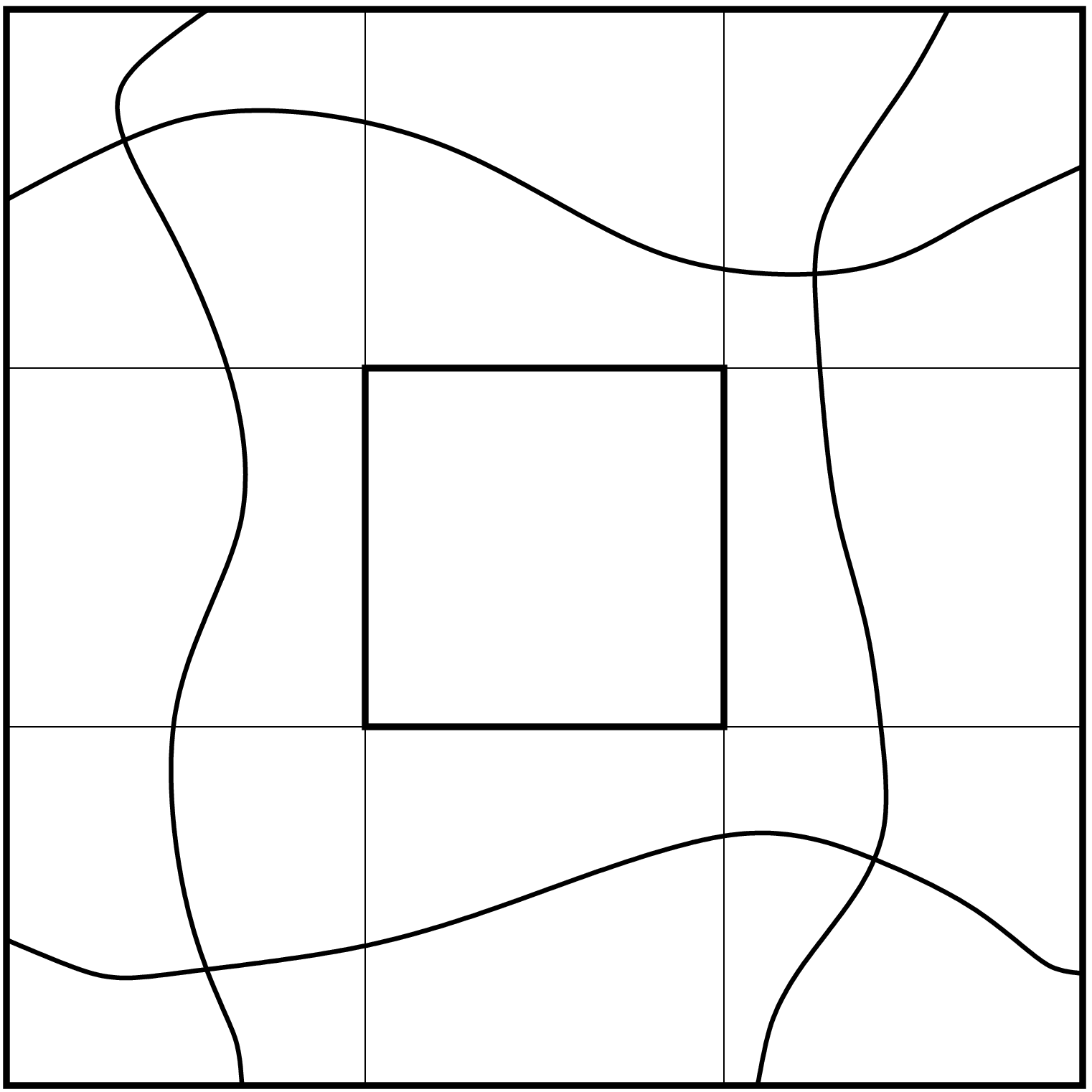}
\caption{Four crossings of rectangles forming a circuit inside an annulus.}
\label{fig3}
\end{center}
\end{figure}

Having established that $\lim_{n \to \infty} h_n = 1$, both the existence and the
uniqueness of an unbounded component with probability one follow from standard pasting
arguments. To prove existence one can use the event depicted in Fig.~\ref{fig2} (and
its rotation by 90 degrees) to couple the complement of the soup to a one-dependent
bond percolation process on a square grid with parameter $p_n \to 1$ as $n \to \infty$
(see Fig.~\ref{fig2}). Thus, choosing $n$ sufficiently large implies percolation in the
bond percolation process and, by the coupling, existence of an unbounded component in the
complement of the soup. Now note that the event ${\cal A}_n$ that the complement of the
soup contains a circuit inside $[-3^{n+1},3^{n+1}]^2 \setminus [-3^n,3^n]^2$ surrounding
$[-3^n,3^n]^2$ has probability bounded below by $h_n^4$ (see Fig.~\ref{fig3}). Hence,
${\cal A}_n$ occurs for infinitely many $n$, which implies the uniqueness of the unbounded
component. \fbox{} \\

\subsection{Proofs Concerning Mandelbrot's Fractal Percolation Model}

\noindent{\bf Sketch of the Proof of Theorem~\ref{thm-mandelbrot}} \,
The fact that $0<\hat{p}_c(N,d)<1$ follows from the same arguments that
show that $0<p_c(N,d)<1$, where $p_c(N,d)$ is the critical probability
defined in terms of crossings of a cube (see Section~\ref{sec-mandelbrot}).

Showing that $\phi_A(\hat{p}_c)>0$ follows the strategy of the proof of
Lemma~\ref{main-lemma}. In fact, the proof is even easier in this case,
since the strategy is particularly well-suited for Mandelbrot's fractal
percolation model. The reason for this lies in the geometry of the fractal
construction.

Using the fact that $\phi_A(p)=0$ for $p<\hat{p}_c$, the proof that
the limiting retained set ${\mathcal C}_N$ is totally disconnected if
$p<\hat{p}_c$ is essentially the same as the proof of Lemma~\ref{lemma2},
to which we refer the reader.
This also shows that $\hat{p}_c(N,d)=\tilde{p}_c(N,d)$. \fbox{} \\

\noindent{\bf Proof of Corollary~\ref{cor-mandelbrot}} \,
Since full space fractal percolation is obtained by tiling ${\mathbb R}^d$
with independent copies of fractal percolation in $[0,1]^d$, it immediately
follows from Theorem~\ref{thm-mandelbrot} that the system is totally
disconnected when $p<\hat{p}_c=\hat{p}_c(N,d)$.

Let us now show that, when $p \geq \hat{p}_c$, for any $M>0$, the system
contains a connected component with diameter larger than $M$ with probability
one. Consider the event that the unit cube $[0,1]^d$ contains a connected
component of diameter larger than $\varepsilon$, and let $\pi(p,\varepsilon)$
denote its probability. Fix an $0<\varepsilon_0<1$ and observe that
by Theorem~\ref{thm-mandelbrot} $\pi(p,\varepsilon_0)>0$. Let $k_0$ be
the smallest integer such that $\varepsilon_0 N^{k_0} > M$.

Consider the full space fractal percolation process obtained by conditioning
on total retention in the unit cube $[0,1]^d$ of the first $k_0$ iterations.
Since the limiting retained set of this process is clearly stochastically
larger than the limiting retained set of the original one, the probability
that the conditioned process contains a connected component of diameter larger
than $\varepsilon_0$ in the unit cube $[0,1]^d$ is at least $\pi(p,\varepsilon_0)>0$.
However, by scaling, this probability is the same as the probability that the
original process contains a connected component of diameter larger than
$N^{k_0} \varepsilon_0 > M$ in the cube $N^{k_0} [0,1]^d$. We can now use
translation invariance to conclude that a connected component with diameter
larger than $M$ is present in the original full space fractal percolation
system with probability one. \fbox{} \\

The next proof relies on ideas developed in~\cite{bc}, and is similar to the
proof of Theorem 1.1 there. For this reason, we present here only a sketch of
the proof, referring the interested reader to~\cite{bc} for more details.

\bigskip

\begin{figure}[!ht]
\begin{center}
\includegraphics[width=5cm]{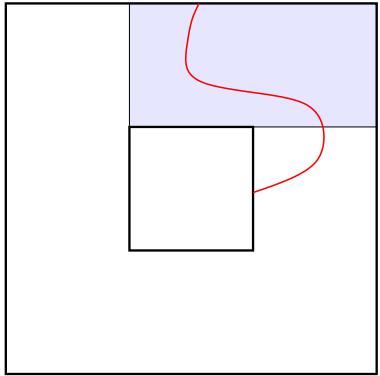}
\caption{The shaded rectangle is crossed in the ``easy'' direction.}
\label{fig4}
\end{center}
\end{figure}

\noindent{\bf Sketch of the Proof of Theorem~\ref{thm-equivalence}} \,
To prove the first statement, it suffices to show that if the shell $A$ in
the definition of $\hat{p}_c$ is crossed with positive probability, then
the unit square $[0,1]^2$ is also crossed with positive probability.
It is easy to see (Fig.~\ref{fig4}) that any crossing of $A$ must cross
either a square of side length $1/3$ or a rectangle of side lengths $1/3$
and $2/3$ in the ``easy'' direction.
Using Theorem \ref{thm-mandelbrot}, this implies that when $p \geq \hat{p}_c$,
there is positive probability of having a crossing of the rectangle
$[0,1/3]\times [0,2/3]$ in the horizontal direction. According to results
from~\cite{dm} (see Lemma~5.1 there) this implies that there is a crossing
of $[0,1] \times [0,2/3]$ with positive probability, which implies the first
statement of the theorem.

Since the limiting retained set ${\mathcal C}_N$ is totally disconnected with
probability one when $p<\hat{p}_c(N,d)$, it is immediate that
$\hat{p}_c(N,d) \leq p_c(N,d)$ for all $N$. Assume that
$p_c(N,d)>\hat{p}_c(N,d)$ for all $N$. Then, for each $N$, we can choose
$p_0=p_0(N,d)$ such that $\hat{p}_c(N,d)<p_0(N,d)<p_c(N,d)$. We will show
that this leads to a contradiction for $N$ large enough.

Consider Mandelbrot's fractal percolation process in $[0,1]^d$ with retention probability $p_0$,
and denote by ${\cal A}_k$ the event that there is complete retention up to the
$k$-th iteration, i.e., ${\mathcal C}^k_N = [0,1]^d$. Let ${\mathbb L}^d$ be the
$d$-dimensional lattice with vertex set ${\mathbb Z}^d$ and with edge set given
by the adjacency relation: $(x_1, \ldots, x_d)=x \sim y=(y_1, \ldots, y_d)$ if
and only if $x \neq y$, $|x_i-y_i| \leq 1$ for all $i$ and $x_i=y_i$ for at least
one value of $i$.

Conditioned on ${\cal A}_{k-1}$, we can couple level $k$ of the fractal percolation
process to a diminishment percolation process (see~\cite{ag,grimmett-book}) on
${\mathbb L}^d$ with the following diminishment rule: for a vertex $x \in {\mathbb Z}^d$,
if all ${\mathbb L}^d$-neighbors of $x$ are closed, except possibly two nearest neighbors
in ${\mathbb Z}^d$, we make $x$ closed, regardless of its state before the diminishment.
Note that this has the effect of ``diminishing" the percolation configuration by changing
the state of some vertices from open to closed. The diminishment is \emph{essential} in
the language of~\cite{ag} (see also~\cite{bc}). The coupling is exactly the same as in
the proof of Theorem~1.1 of~\cite{bc} (although the diminishment rule is different),
therefore the interested reader can check the details in~\cite{bc}.

Let $\psi_k(p_0)$ denote the probability that there is
an open crossing between the inner and outer boundaries of $N^k A$ in the diminishment
percolation process with initial density $p_0$ of open vertices. A feature of the
coupling is that the diminishment percolation process dominates the fractal percolation
process in the sense that $\phi_A(p_0) \leq \psi_k(p_0)$ (see~\cite{bc}).

From~\cite{FG1,FG2} we know that, for all $d \geq 2$, $p_c(N,d) \to p'_c(d)$ as
$N \to \infty$, where $p'_c(d)$ is the critical value for Bernoulli site percolation
on ${\mathbb L}^d$. Furthermore, for the diminished percolation model, the critical value $p_c''(d)$ satisfies
$p_c''(d)>p'_c(d)$ (which follows from the diminishment being essential,
see~\cite{ag,grimmett-book} for more details on enhancement
and diminishment percolation).

This implies that, for fixed $d$ and $N$ sufficiently large,
the diminishment percolation process with initial density $p_0=p_0(N,d)<p_c(N,d)<p_c''(d)$.
Observe that it is not subcritical in the sense that
$p_0<p'_c(d)$, rather it is the {\em diminished} percolation process that is
subcritical. Therefore, for $N$ sufficiently large,
$\lim_{k \to \infty} \psi_k(p_0) = 0$, which implies that $\phi_A(p_0)=0$.
However, since $p_0>\hat{p}_c(N,d)$, $\phi_A(\hat{p}_c(N,d))>0$ by Theorem~\ref{thm-mandelbrot},
and crossing probabilities are increasing in $p$, this leads to a contradiction. \fbox{} \\

\noindent{\bf Remark} \, The result from \cite{dm} that is used for the above proof
when $d=2$ uses planar arguments and therefore cannot be readily generalized to $d\geq 3$.

\bigskip

\noindent {\bf Acknowledgements} \, The second author thanks the Department
of Mathematics of Chalmers University of Technology for the hospitality
during a visit in 2008 when this work was started. He also thanks the
Netherlands Organization for Scientific Research (NWO) for making the
visit possible and for financial support through a VENI grant while
part of this work was carried out. The first author thanks the Department
of Mathematics of the Vrije Universiteit Amsterdam for the hospitality
during a subsequent visit. Both authors thank Institut Mittag-Leffler
for the hospitality during the Spring 2009 program on Discrete Probability
when this paper was completed. They also thank Ronald Meester for
suggesting to look at the multiscale Boolean model studied in~\cite{mr}
and for giving many useful comments on the manuscript. Finally they thank
an anonymous referee who provided many good comments and useful suggestions.

\end{document}